# Vertex degree distribution and arc endpoints degree distribution of graphs with a linear rule of preferential attachment and Pennock graphs

V N Zadorozhnyi, E B Yudin

Omsk state technical university, pr. Mira 11, 644050, Omsk, Russia

E-mail: zwn2015@yandex.ru, udinev@asoiu.com

**Abstract.** The article deals with two classes of growing random graphs—graphs following the preferential attachment rule with a linear weight function (*L*-graphs) and hybrid Pennock graphs. We determine the exact final vertex degree distribution and the exact final two-dimensional arcs/edges degree distributions of graphs under consideration. The study prove that each hybrid Pennock graph is isomorphic to a certain *L* graph, and that the converse does not hold, since there are no Pennock graphs isomorphic to *L* graphs with negative displacements in the expression for the linear weight function. A formula is derived that makes it possible to determine the weight functions, which are used to generate graphs with the required asymptotic power-law vertex degree distribution. The reliability of the obtained results is confirmed by calculations using accurate numerical methods and simulation modeling, i.e. by direct generating of the graphs. The practical value of the results is demonstrated by an example of their effective application for accurate calibration of a growing graph that simulates a network of Internet at the level of autonomous systems.

## 1. Introduction

The preferential attachment graphs spread widely as mathematical models of large networks after the publication of the article [1], in which a growing graph was proposed, later called Barabasi-Albert graph (BA graph). The BA graph, being an adequate model of many growing networks, made it possible to explain how such structural properties arise as the asymptotic power-law vertex degree distribution (VDD), the error and attack tolerance, properties of small-world networks (six degrees of separation), etc. [2].

In this article we consider both the growing graphs of *L* class – graphs with a linear weight function, which is written in the form $f(k) = k + s$, where the displacement *s* is a constant, either in the form $f(k) = s = 1$, and the graphs of *P* class, i.e. Pennock graphs [3].

*Graphs of L class* are grown from a small seed-graph by adding increments — new vertices with a given number *m* of outgoing arcs. Each free end of the every added increment arc attaches to *i*-th vertex of the graph with probability $p_i$, proportional to the weight $f(k_i)$ of this vertex, defined by its vertex degree $k_i$:

$$p_i = \frac{f(k_i)}{\sum_{j=1}^{N} f(k_j)} = \frac{k+s}{\sum_{j=1}^{N}(k_j + s)}, \qquad i = 1,...,N, \qquad (1)$$

where $k = k_i$ and $N$ is a current number of vertices in a graph. From the first equality in (1) it is seen that all weight functions $f_c(\cdot) = cf(\cdot)$, which differ from the function $f(\cdot)$ by the multiplicative coefficient $c > 0$, are equivalent and determine the same probability $p_i$. It follows that the «special» weight function $f(k) = 1$ of $L$ class graphs is a typical weight function $f(k) = k + s$, which is equivalent to the weight function $f_c(k) = f_{1/s}(k) = k/s + 1$ obtained in the limit as $s \to \infty$. Thus, every graph of $L$ class is given by two parameters – degree $m$ of the vertex in the increment and displacement $s$. BA graph is determined by the weight function $f(k) = k$, i.e. it is such a graph of $L$ class whose weight function has zero displacement. In general case, the displacement $s$ can take negative values.

*Every graph of P class* – Pennock graph – is also given by two parameters – degree $m$ of the vertex in the increment and the probability $a$, which determines a fraction of the weight function $f_1(k) = 1$ in its «mixture» with the weight function $f_2(k) = k$. When linking to the graph according to rule (1) each arc of an increment uses a weight function $f_1(k) = 1$ with probability $a$, and a weight function $f_2(k) = k$ with probability $(1 - a)$.

The aim of the study is to detect such graphs in the $L$ and $P$ classes whose final (obtained as $N \to \infty$) distributions $Q_k$ of vertex degree $k$ have power-law asymptotics, i.e. for which $Q_k \sim ck^{-\alpha}$ as $k \to \infty$ (the notation $c$ of an indefinite numerical coefficient used by us in different contexts has different meanings). At the same time the exponent $\alpha$ must belong to the interval $2 < \alpha \leq 3$, as in many real growing networks, for example, in the Internet at the level of autonomous systems. The probability $Q_k$ of degree $k$ here is understood to be the probability, that a vertex chosen randomly (with equal probability) in an infinite graph has a degree $k$. The average degree of vertices of the growing graphs under consideration, as follows from their construction, is equal to value $\langle k \rangle = 2m$, this is the number of edge ends per one vertex of an infinite graph.

The hybrid Pennock graphs are interesting because as $a \to 1$, as it follows from the algorithm of their construction, they converge to a known graph with geometric VDD $Q_k$, as $a \to 0$ – to BA graph, i.e. to the graph with asymptotically-power VDD $Q_k \sim ck^{-3}$. By choosing the value of the parameter $a$ of the hybrid graph between 0 and 1, we can obtain graphs with «intermediate» types of VDD [3, 4].

## 2. Final VDD in *P*- and *L*-graphs

Let us apply the method of derivation of exact recurrence relations to the hybrid Pennock graph with parameter $a$, described in [5], adapting it to the features of the hybrid graph growing algorithm. Namely, we divide the set of vertices of the hybrid graph into layers (subsets) $A_k$ so that each layer $A_k$ consists of vertices with a degree $k$. When adding the next increment to the graph, the probability $p_i$ that a new arc attaches to a vertex $i$ (having a degree $k$) is

$$p_i = p(k) = a \frac{f_1(k)}{\sum_{j=1}^{N} f_1(k_j)} + (1-a) \frac{f_2(k)}{\sum_{j=1}^{N} f_2(k_j)} \sim a \frac{f_1(k)}{\langle f_1 \rangle N} + (1-a) \frac{f_2(k)}{\langle f_2 \rangle N}, \qquad (2)$$

where $\langle f_1 \rangle$ is an average value of the weight $f_1$; $\langle f_2 \rangle$ is an average value of the weight $f_2$.

As $N$ grows, the average number of vertices of the $A_k$ layer selected by $m$ arcs of an increment to attach to converges (with relative error of zero) to the value

$$\Delta_k^- = mp(k)|A_k| = mp(k)Q_k N = am \frac{f_1(k)}{\langle f_1 \rangle} Q_k + (1-a)m \frac{f_2(k)}{\langle f_2 \rangle} Q_k.$$

All the vertices chosen in the $A_k$ layer move to the next layer $A_{k+1}$, since the addition of new arcs to them increases their degree by one. We neglect the probability of attaching more than one arc of an increment to a vertex $i$, since it is of infinitesimal order $o(p(k))$ for the functions under consideration $f_1(k) = 1$ and $f_2(k) = k$ [5]. Due to the departure of the vertices selected in the $A_k$ the average number of vertices in this layer decreases by $\Delta_k^-$. Similarly, we determine the average number of vertices selected by increment in the layer $A_{k-1}$, which is

$$\Delta_k^+ = mp(k-1)Q_{k-1}N = am\frac{f_1(k-1)}{\langle f_1 \rangle}Q_{k-1} + (1-a)m\frac{f_2(k-1)}{\langle f_2 \rangle}Q_{k-1}.$$

These vertices move to the $A_k$ layer thereby increasing the average number of vertices in it.
In the limit at $t \to \infty$ ($N \to \infty$) the probabilities $Q_k = |A_k|/N = N_k/N$ do not change after adding an increment to the graph (a new vertex with $m$ arcs), as they are stationary probabilities, i.e. the following two equalities (balance equations) are valid:

$$\frac{N_k}{N} = \frac{N_k - \Delta_k^- + \Delta_k^+}{N+1} \text{ (at } k \geq m+1\text{)}, \qquad \frac{N_m}{N} = \frac{N_m - \Delta_k^- + 1}{N+1} \text{ (at } k = m\text{)}.$$

The last equation is the equation for the layer $A_m$, it differs from the previous one because no vertices come into $A_m$ layer from the previous $A_{m-1}$ layer (there is no $A_{m-1}$ layer), but at the same time, the vertex of the graph increment is added to the $A_m$ layer.

We substitute the above expressions for $\Delta_k^-$ and $\Delta_k^+$ in the first balance equation:

$$\frac{N_k}{N} = \frac{1}{N+1}\left[N_k - \left(am\frac{f_1(k)}{\langle f_1 \rangle}Q_k + (1-a)m\frac{f_2(k)}{\langle f_2 \rangle}Q_k\right) + am\frac{f_1(k-1)}{\langle f_1 \rangle}Q_{k-1} + (1-a)m\frac{f_2(k-1)}{\langle f_2 \rangle}Q_{k-1}\right].$$

Multiplying both parts of the resulting equation by $N(N+1)$, giving similar terms and again dividing both parts of the equation by $N$, we rewrite it, taking into account that $N_k/N \sim Q_k$, in the form

$$Q_k = am\frac{f_1(k-1)}{\langle f_1 \rangle}Q_{k-1} + (1-a)m\frac{f_2(k-1)}{\langle f_2 \rangle}Q_{k-1} - am\frac{f_1(k)}{\langle f_1 \rangle}Q_k - (1-a)m\frac{f_2(k)}{\langle f_2 \rangle}Q_k.$$

Expressing now $Q_k$ via $Q_{k-1}$, we obtain the recurrent formula

$$Q_k = \frac{am\frac{f_1(k-1)}{\langle f_1 \rangle} + (1-a)m\frac{f_2(k-1)}{\langle f_2 \rangle}}{1 + am\frac{f_1(k)}{\langle f_1 \rangle} + (1-a)m\frac{f_2(k)}{\langle f_2 \rangle}}Q_{k-1}, \qquad k \geq m+1. \qquad (3)$$

Similarly, from the second balance equation for the layer $A_m$ we find:

$$Q_m = \frac{1}{1 + am\frac{f_1(m)}{\langle f_1 \rangle} + (1-a)m\frac{f_2(m)}{\langle f_2 \rangle}}. \qquad (4)$$

Formulas (3), (4) are valid not only for hybrid Pennock graphs, since their derivation is not connected with the type of mixed weight functions. For the Pennock graphs, substituting the expressions for the used weight functions $f_1(k) = 1$, $f_2(k) = k$ and their mean values $\langle f_1 \rangle = \langle 1 \rangle = 1$, $\langle f_2 \rangle = \langle k \rangle = 2m$ in (3), (4), we obtain the following recurrence formulas that determine VDD $Q_k$ exactly:

$$Q_m = \frac{2}{2 + am + m}, \qquad (5)$$

$$Q_k = \frac{2am+(1-a)(k-1)}{2+2am+(1-a)k} Q_{k-1}, \qquad k \geq m+1. \qquad (6)$$

The validity of formulas (5), (6) is confirmed by simulation modeling (direct cultivation of P-graphs) and by a comparison of the resulting VDD with the VDD calculated by formulas (5), (6).

Now let us consider the graphs of L class. Applying the derived in [6] recurrent formulas, which are common for preferential attachment graphs,

$$Q_m = \frac{\langle f \rangle}{\langle f \rangle + mf(m)}, \qquad (7)$$

$$Q_k = \frac{mf(k-1)}{\langle f \rangle + mf(k)} Q_{k-1}, \qquad k \geq m+1, \qquad (8)$$

i.e. substituting the weight function $f(k) = k + s$ and the average weight $\langle f \rangle = \langle k + s \rangle = \langle k \rangle + s = 2m + s$ in these formulas, we obtain:

$$Q_m = \frac{2m+s}{2m+s+m(m+s)}, \qquad (9)$$

$$Q_k = \frac{m(k-1+s)}{2m+s+m(k+s)} Q_{k-1}, \quad k \geq m+1. \qquad (10)$$

For a linear function $f(k) = 1$ the formulas (7), (8) take the form

$$Q_m = \frac{1}{1+m}, \qquad (11)$$

$$Q_k = \frac{m}{1+m} Q_{k-1}, \quad k \geq m+1, \qquad (12)$$

i.e. define a decreasing geometric progression with the initial term $1/(1+m)$ and common ratio $m/(1+m)$. Formulas (11), (12) are obtained from formulas (9), (10) as $s \to \infty$.

Now we can prove the following statement connecting two classes of random graphs under consideration.

T h e o r e m 1. Every hybrid graph of P class is equivalent in the final VDD to a certain graph of L class.

P r o o f. The graph of P class and its VDD is uniquely determined by the values of the parameters $m$ and $a$. In accordance with any selected P- graph with the parameter $a < 1$ we put a L-graph with the same $m$ and with the parameter

$$s = \frac{2am}{1-a}. \qquad (13)$$

Substituting these $m$ and $s$ into formulas (9), (10) and simplifying them, we obtain formulas that coincide, respectively, with formulas (5), (6). Hence, the VDD of the L-graph with such $m$, $s$ coincides with the VDD of the selected P-graph. Let the hybrid graph with the parameter $a = 1$ correlate with the L-graph having a weight function $f_1(k) = 1$. In this case, the VDD of the L-graph (11), (12) also coincides with the VDD (5), (6) of the selected hybrid graph for which $a = 1$. The theorem is proved.

The reverse is not true, since no P-graphs correspond to L-graphs with the parameter $s < 0$. Therefore, to search for graphs with asymptotic power-law VDD in the P and L classes it is necessary and sufficient to consider graphs of L class.

## 3. Exact formulas for arc endpoints degree distribution and edge endpoints degree distribution in *P*-graphs

An arc endpoints degree is an ordered pair of numbers $(l, k)$, where $l$ is a degree of a vertex from which the arc emanates, $k$ is a degree of a vertex into which the arc enters. The probability $Q_{l,k}$ is defined as the probability that a randomly (with equal probability) selected arc of an infinite graph will have a endpoints degree $(l, k)$. The set of all probabilities $Q_{l,k}$ ($l \geq 1, k \geq 1$), i.e. the final two-dimensional arc endpoints degree distribution of *P*-graph can be found by compiling and solving the balance equations for the arc endpoints degree for the graphs with the nonlinear rule of preferential attachment (NPA rule) [5]. The probability $p(k)$ of attaching an arc of an increment to a vertex of degree $k$ should be determined by the formula of the total probability (2). This entails a change in the formula for the probability $P_k$ of attaching an arc of an increment to a layer $A_k$: in *P*-graphs $P_k = [a + k(1-a)/(2m)]Q_k$. Omitting routine calculations in connection with the limited volume of the article, we present the resulting arc endpoints degree distribution for the *P*-graph, given by parameters $m$, $a$

$$Q_{l,k} = \begin{cases} 0 & \text{при } l < m \text{ или } k \leq m, \\ \dfrac{[2am + (1-a)(k-1)] \cdot [Q_{k-1} + mQ_{m,k-1}]}{2m + (3a+1)m^2 + m(1-a)k} & \text{при } l = m, k \geq m+1, \quad (14) \\ \dfrac{[2am + (1-a)(k-1)]Q_{l,k-1} + [2am + (1-a)(l-1)]Q_{l-1,k}}{2 + 4am + (1-a)(k+l)} & \text{при } l \geq m+1, k \geq m+1. \end{cases}$$

The calculation of the array $\mathbf{Q} = \| Q_{l,k} \|$ starts with zero filling of rows with numbers $l < m$ (if $m > 1$) and columns with numbers $k \leq m$ according to the first row of the solution (14). Then, using the second solution formula (14) the array row $l = m$ is calculated starting from the $(m+1)$-th element, from left to right, and the third solution formula (14) is used to calculate the next rows of the array, also from the $(m+1)$-th element.

Solution (14) is verified and validated by simulation. It is easy to check the solution (14) in particular cases for $a = 0$ or $a = 1$. In these cases, the hybrid *P*-graph degenerates into an ordinary graph of preferential attachment: for $a = 1$ – into a graph with the weight function $f(k) = f_1(k) = 1$, for $a = 0$ – into a graph with the weight function $f(k) = f_2(k) = k$. It is easy to verify that in these cases the solution (14) coincides with the solution found in [5] for graphs with NPA rule:

$$Q_{l,k} = \begin{cases} 0, & l < m \text{ или } k \leq m, \\ \dfrac{f_m Q_m}{\langle f \rangle + m(f_m + f_{m+1})}, & l = m, k = m+1, \\ \dfrac{f_{k-1}(Q_{k-1} + mQ_{m,k-1})}{\langle f \rangle + m(f_m + f_k)}, & l = m, k \geq m+2, \quad (15) \\ \dfrac{f_{l-1}Q_{l-1,k} + f_{k-1}Q_{l,k-1}}{\langle f \rangle / m + f_l + f_k}, & l \geq m+1, k \geq m+1, \end{cases}$$

where the function $f(k)$ of the integer variable $k$ is denoted by $f_k$.

If all arcs in a directed graph are replaced by edges, then we get an undirected graph. The edge endpoints degree array $\Theta$ is obtained as a result of the following transformation of the arc endpoints degree array $\mathbf{Q}$:

$$\Theta = (\mathbf{Q} + \mathbf{Q}^T)/2. \quad (16)$$

Here, the symbol T denotes the transposition of the array **Q**. The probability $\Theta_{l,k}$ (an element of the array $\Theta$) is a probability that a randomly selected edge of an infinite graph traversed in a random direction, leads from a vertex with degree $l$ to a vertex with degree $k$. The array $\Theta$ is symmetric with respect to the main diagonal.

### 4. Exact formulas for arc endpoints degree distribution and edge endpoints degree distribution in *L*-graphs.

The exact formulas for arc endpoints degree distribution of graphs with the linear rule of preferential attachment are derived as a special case of the general formula for VDD (15) by substituting a linear weight function $f(k) = k + s$ and, respectively, the average weight of the vertices $\langle f \rangle = 2m + s$:

$$Q_{l,k} = \begin{cases} 0, & l < m \text{ или } k \leq m, \\ \dfrac{(m+s)Q_m}{2m+s+m(2m+2s+1)}, & l = m, k = m+1, \\ \dfrac{(k-1+s)(Q_{k-1}+mQ_{m,k-1})}{2m+s+m(m+2s+k)}, & l = m, k \geq m+2, \\ \dfrac{(l-1+s)Q_{l-1,k}+(k-1+s)Q_{l,k-1}}{(2m+s)/m+l+2s+k}, & l \geq m+1, k \geq m+1. \end{cases} \quad (17)$$

T h e o r e m 2. Every hybrid graph of *P* class, equivalent in the final VDD to a definite graph of *L* class, is also equivalent to this graph in the final arc/edge endpoints degree distribution.

P r o o f. According to Theorem 1 any *P*-graph with parameters $m$, $a$ is equivalent in the VDD $Q_k$ to *L*-graph with the same $m$ and the parameter

$$s = 2am/(1-a). \quad (18)$$

Then, let us verify that these *P*- and *L*-graphs are also equivalent with respect to arc endpoints degree distribution $Q_{l,k}$ and edge endpoints degree distribution $\Theta_{l,k}$. Taking into account formula (16) it is sufficient to check the equivalence of these graphs with respect to arc endpoints degree distribution. To do this, we substitute the value $s$ (18) into the general formula of arc endpoints degree distribution for the *L*-graph (17). Then, after equivalent transformations, formula (17) reduces to (14). Consequently, the compared *P*- and *L*-graphs are equivalent in arc endpoints degree distribution $Q_{l,k}$. The theorem is proved.

Figure 1 shows the examples of edge endpoints degree distributions for two equivalent *P*- and *L*-graphs.

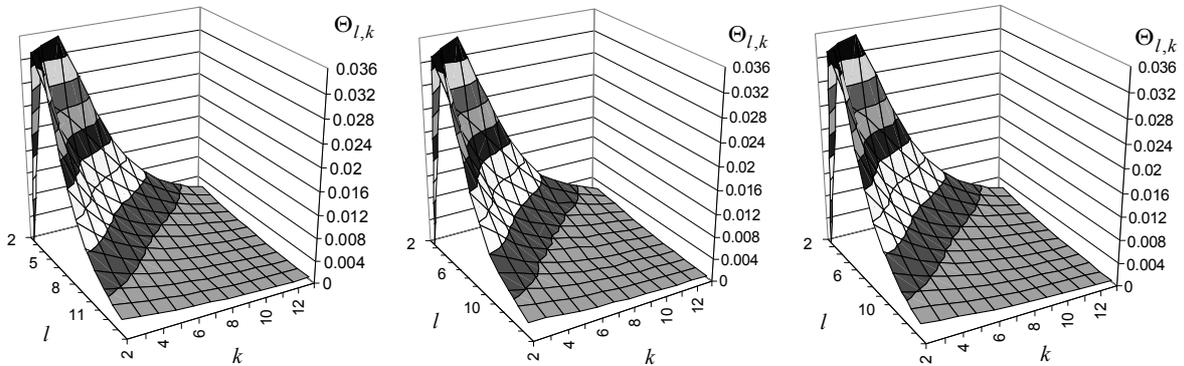

**Figure 1**. Left: edge endpoints degree distribution (16), (14) $\Theta_{l,k}$ of Pennock graph with parameters $m = 2$, $a = 0.75$; center: edge endpoints degree distribution of the same graph, calculated by its simulation; right: VDD of equivalent *L*-graph with parameters $m = 2$, $s = 2am/(1-a) = 12$

The graphs in figure 1 clearly confirm the validity of the formulas obtained for the edge endpoints degree distribution (and, therefore, for the formulas used for VDD and arc endpoints degree distribution) of the *P*- and *L*-graphs, as well as the validity of the theorems proved in this article

## 5. Theorem which states the isomorphism of Pennock graph models and L-models with displacements $s \geq 0$

The exact coincidence of the distributions $Q_k$ and $Q_{lk}$ in the *P*-graph with parameters $(m, a)$, on the one hand, and in the *L*-graph with parameters $(m, s = 2am/(1 - a))$ on the other hand, leads to the assumption of a complete isomorphism of the compared $P(m, a)$ and $L(m, s)$ network models. This assumption is confirmed below. Its validity is proved by the following theorem.

T h e o r e m  3. Every *P*-graph is equivalent to a certain *L*-graph.

P r o o f . Let the *P*-graph be given by the parameters $m$ (a number of arcs in a graph increment) and $a$ (a probability of using the weight function $f_1(k) = 1$ when selecting a vertex) and let the *L*-graph have the same parameter $m$ and the weight function $f(k) = k + s$, where $s = 2am/(1-a)$. Then the rule (2) of choosing a vertex to attach to with a new arc of a *P*-graph is equivalent to the rule (1) for selecting a vertex with a new arc of an *L*-graph.

Indeed, simplifying the rule (2) of choosing a vertex in a *P*-graph, we find that

$$p_i = a \frac{f_1(k)}{\sum_{j=1}^{N} f_1(k_j)} + (1-a) \frac{f_2(k)}{\sum_{j=1}^{N} f_2(k_j)} = a \frac{1}{\sum_{j=1}^{N} 1} + (1-a) \frac{k}{\sum_{j=1}^{N} k_j} = a \frac{1}{N} + (1-a) \frac{k}{2R} =$$

$$= a \frac{1}{N} + (1-a) \frac{k}{2mN} = \frac{2am + k(1-a)}{2mN}, \qquad i = 1,...,N, \qquad (19)$$

where $R = mN$ is a number of edges in the graph. In the calculation (19) for simplicity we assume that the number of edges in the used seed-graph is $m$ times greater than the number of vertices in it (as in any graph increment). Now calculating the probability $p_i$ for the corresponding *L*-graph according to the rule (1), we obtain

$$p_i = \frac{k+s}{\sum_{j=1}^{N}(k_j + s)} = \frac{k + \frac{2am}{1-a}}{\sum_{j=1}^{N} \left(k_j + \frac{2am}{1-a}\right)} = \frac{k(1-a) + 2am}{\sum_{j=1}^{N} k_j(1-a) + \sum_{j=1}^{N} 2am} = \frac{k(1-a) + 2am}{(1-a)\sum_{j=1}^{N} k_j + 2amN} =$$

$$= \frac{k(1-a) + 2am}{(1-a)2R + 2amN} = \frac{k(1-a) + 2am}{2(1-a)mN + 2amN} = \frac{k(1-a) + 2am}{2mN}, \qquad i = 1,...,N. \qquad (20)$$

Consequently, the rules of growing the compared *P*- and *L*- graphs are equivalent (probabilities (19) and (20) of attaching to any *i* vertex coincide).

The theorem is proved.

Hence, *P*-graphs are not a new class of graphs – the set of *P*-graphs coincides with the subset of *L*-graphs, defined by the condition $s \geq 0$. Any *P*-graph is equivalent to the corresponding *L*-graph in all probabilistic and probabilistic-temporal characteristics, in particular, in the characteristics of the transition process. Of course, the characteristics of transient processes are equivalent for the same seed-graphs, from which the compared *P*- and *L*-graphs are grown.

## 6. Derivation of an approximate expression in closed form for VDD $Q_k$ in the *L* class

For *L*-graphs we find the expressions of VDD $Q_k$ in closed form, which allow to analyzing the asymptotic behavior of these VDD.

In the case of a linear weight function $f(k) = 1$ recursion (11), (12) implies an exact explicit formula of VDD, this is a formula for a general term of the geometric progression defined by the equalities (11), (12):

$$Q_k = \frac{1}{1+m}\left(\frac{m}{1+m}\right)^{k-m}, \quad k \geq m. \tag{21}$$

In the case of a linear weight function $f(k) = k + s$, it is also easy to obtain an exact expression for $Q_k$ in closed form, since the corresponding recurrent equations (9), (10) are linear with respect to the desired function $Q_k = Q(k)$. Solving equations (9), (10), we obtain the following exact expression:

$$Q_k = \frac{(2m+s)\Gamma(m+s+2+s/m)\Gamma(k+s)}{m\Gamma(m+s)\Gamma(k+s+3+s/m)}. \tag{22}$$

As $s \to \infty$ we find that $Q_k$ (22) converges to $Q_k$ (12). The solution (22) of equations (9), (10) is easily proved by the method of mathematical induction with a help of the known identity $\Gamma(z+1) = z\Gamma(z)$. However, despite the closed form of the resulting expression (22), it is inconvenient for the analysis of the asymptotic of VDD $Q_k$ as $k \to \infty$, obtained at different $s$.

Therefore, we proceed to the derivation of an approximate expression for the VDD, which we find using the mean-field theory [1, 3, 4]. The vertices added to the graph will be numbered by the moments $i$ of their arrival. The number $N$ of vertices in the graph will always be equal to the current time $t$ of the next increment. With this in mind, we perform the following steps prescribed by the mean-field theory, which have a fairly simple probabilistic interpretation.

*Step* 1). Based on (1) we write down the differential equation for the time variation of an average degree $k_i$ of a vertex $i$ due to attaching of new graph increments:

$$\frac{dk_i}{dt} = \frac{mf(k_i)}{\langle f \rangle t} = \frac{m(k_i + s)}{\langle f \rangle t}.$$

Since here $\langle f \rangle = \langle k + s \rangle = 2m + s$, then

$$\frac{dk_i}{dt} = \frac{m(k_i + s)}{(2m+s)t}.$$

Solving this equation with separable variables, taking into account the initial condition $k_i(i) = m$, we obtain

$$k_i = (m+s)\left(\frac{t}{i}\right)^{\frac{m}{2m+s}} - s. \tag{23}$$

*Step* 2). Using (23), we find the number $i = l$ of such a vertex whose degree is equal to $k$:

$$\frac{l}{t} = \left(\frac{k+s}{m+s}\right)^{-\frac{2m+s}{m}}.$$

*Step* 3). Hence, we obtain an estimate $\hat{F}(k)$ of a distribution function (d.f.) $F(k)$ of a degree $k$ for a randomly chosen vertex in the graph:

$$\hat{F}(k) = \frac{t-l}{t} = 1 - \left(\frac{k+s}{m+s}\right)^{-\frac{2m+s}{m}}. \tag{24}$$

*Step* 4). An estimate for the VDD is found as a derivative of d.f. (24):

$$\hat{Q}_k = \frac{2m+s}{m}(m+s)^{\frac{2m+s}{m}}(k+s)^{-\frac{3m+s}{m}}, \tag{25}$$

where $\hat{Q}_k$ – is an approximate (with asymptotically exact decay rate) estimate for $Q_k$.

## 7. Analysis of asymptotic properties of VDD $Q_k$ in *L*-graphs.

The obtained estimate (25) shows that the VDD of the *L*-graph for any finite *s* is asymptotically power-law. This fact simplifies the analysis, since it can be reduced to viewing *s*, multiples of *m* and subsequent interpolation of the obtained results to intermediate *s* values. Let us consider how an asymptotic behavior of the VDD $Q_k$ changes with a change in the parameter *s*.

**1)**. When *s* = 0 from (25) we find that $\hat{Q}_k = 2m^2 k^{-3}$. This value of *s* transforms the weight function $f(k) = k + s$ into the weight function of BA graph and the estimate (25) – into the correct estimate of the asymptotic behavior $Q_k \sim ck^{-3}$ of BA graph. Thus, the value *s* = 0 defines one of the power-law asymptotic distribution $Q_k \sim ck^{-\alpha}$ we are interested in, for which $2 < \alpha \le 3$.

**2)**. When *s* = *m*, 2*m*, 3*m*, … , according to (25), a power-law asymptotic distribution $Q_k \sim ck^{-\alpha}$ with the parameter α = 4, 5, 6, … is realized. Thus, in a semi-infinite range *s* ≥ 0 only one value (the value *s* = 0) provides a power-law asymptotic distribution of the *L*-graph with the value α that belongs to the range $2 < \alpha \le 3$. It is the value α = 3.

Note, that relation (13) establishes one-to-one correspondence between points of the range *s* ≥ 0 and points of the range $0 \le a < 1$. It follows that in *P*-class only one graph satisfies the criterion ($2 < \alpha \le 3$), it is BA graph, realized at *a* = 0.

**3)**. The rest of the range *s* < 0 is limited by the condition of rule feasibility (1): $f(k) = k + s > 0$ for any degree $k \ge m$, it requires that *m* + *s* > 0 must be satisfied, i.e. *s* > –*m*.

When *s* = – 0, as we saw above in paragraph 1), the asymptotic distribution $Q_k \sim ck^{-3}$ takes place.

At the point *s* = – *m* + ε as $\varepsilon \downarrow 0$, according to (25), $Q_k \sim ck^{-\alpha}$ as $\alpha \downarrow 2$. Obviously, the values of the parameter *s* in the interval $-m < s \le 0$ lead to the values α in the interval $-2 < \alpha \le 3$.

Therefore, the purpose of the study is achieved. The desired class of graphs with the asymptotic distribution $Q_k \sim ck^{-\alpha}$, where $2 < \alpha \le 3$, is a class of *L*-graphs with the weight function $f(k) = k + s$, at $m < s \le 0$. If it is necessary to grow a *L*-graph with the given α value, that belongs to $2 < \alpha \le 3$, then the parameter *s* of the weight function $f(k) = k + s$ is determined by the formula that follows from (25)

$$s = (\alpha - 3)m. \qquad (26)$$

*Note* 1. When modeling real networks, the value of *m* is often not integral. The graphs simulating such networks are grown by using stochastic increments: new vertices added to the graph contain a random number *x* arcs with an average value $\langle x \rangle = m$. This does not affect the validity of formulas (25), (26), since the asymptotic VDD is uniquely determined by the weight function and *m* arcs in an increment at the average.

The results of the performed analysis of *P* and *L* graph classes are summarized in table 1.

**Table 1.** Asymptotic VDD $Q_k$ in *P*- and *L*-graphs

| Asymptotic VDD $Q_k$ for $k \to \infty$ | Values of *s* parameter of weight function in *L*-graph | Values of *a* parameter of equivalent *P*-graph |
|---|---|---|
| $Q_k \sim ck^{-\alpha}$, $2 < \alpha < 3$ | $-m < s < 0$ | there are no such *P*-graphs |
| $Q_k \sim ck^{-3}$, $\alpha = 3$ | $s = 0$ | $a = 0$ |
| $Q_k \sim ck^{-\alpha}$, $\alpha > 3$ | $s > 0$ | $0 < a < 1$ |
| exponential | $s \to \infty$ : $f(k) = 1$ | $a = 1$ |

*Note* 2. The value α = 2 in the power-law asymptotic distribution of graphs and networks with preferential attachment is unattainable under any weight functions, since asymptotic distribution

$Q_k \sim ck^{-2}$ leads to an infinite average degree $\langle k \rangle = \sum_{k=m}^{\infty} kQ_k = \infty$, which is impossible, since in graphs (networks) of preferential attachment $\langle k \rangle = 2m$ by construction.

## 8. Example of application

The power-law asymptotic distribution for $Q_k \sim ck^{-\alpha}$ when $2 < \alpha \leq 3$ is typical, for example, for networks that make up the Internet. Paper [6] provides an example of solving the problem of calibrating a graph that simulates a network of autonomous Internet systems according to statistics on this network. The network consists of $N = 22963$ nodes and $E = 48436$ edges. We also have information about the node degree distribution of the network. In total, the table of statistical data contains 1713 lines consisting of two numbers $k$, $n_k$, where $k$ is a degree of a vertex ($k = 1, ..., 1713$), $n_k$ is a number of nodes with a degree $k$. The empirical VDD $Q_k$ is obtained by adding a third column containing numbers $Q_k = n_k / N$ (this VDD is shown in figure 2 by markers).

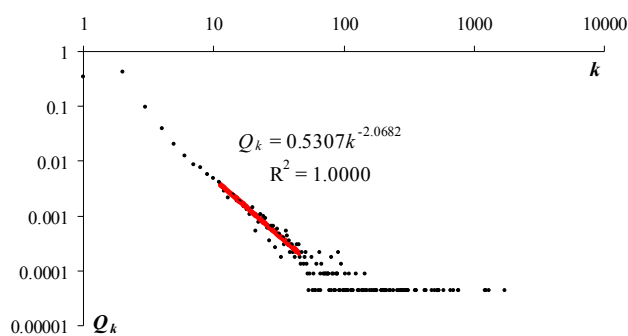

**Figure 2.** Distribution of $Q_k$ on the Internet

The value of the network parameter $m$ is $E/N = 2.1093$. To find the empirical value of the exponent $\alpha$, the power trend line $Q_k \sim ck^{-\alpha}$ (red solid line in figure 2) is constructed, the values $c$ and $\alpha$ are determined by the least squares method (LSM). In this case, the initial part of the empirical VDD is excluded, since it is required to determine the asymptotic characteristic. Empirical probabilities $Q_k$ of large $k$ values are also not used, since there are still few nodes with such degrees in the growing network, and empirical probabilities of such vertex degrees are unreliable (see figure 2). As a result of applying the LSM, we obtain an empirical asymptotic distribution $Q_k \sim ck^{-\alpha}$, in which $\alpha = 2.0682$. By the formula (26) we calculate $s = -1.9655$.

Therefore, the desired asymptotic behavior is realized in a graph with a weight function $f(k) = k - 1.9655$ (with the average number of arcs in the increment of $m = 2.1093$). By defining the weight function that provides the desired asymptotic distribution $Q_k$, we now easily solved the main calibration problem.

In addition, using the methods described in [5–7], by determining the appropriate stochastic increment (when the number of outgoing arcs of the increment is not deterministic, but a random variable), we, as before, solve the problem of implementing a specific initial part of the sequence $Q_k$. As a result, the autonomous system network model is now defined as a graph with preferential attachment, whose weight function when $k \geq 11$ is given by the formula $f(k) = k - 1.9655$, and for $k = 1, ..., 10$ takes values 0, 0, 0, 1.38802, 2.40613, 5.28966, 6.67, 6.71098, 7.79545, 8.10619 respectively. The stochastic increment is formed as a vertex with a random number $x$ of arcs, which can take the values $x = 1, ..., 6$ with probabilities 0.34145, 0.42246, 0.09664, 0.09433, 0.01504, 0.03008 respectively (meanwhile, the average value of $x$ is $m = 2.1093$).

Figure 3 shows the result of validation the performed calibration with a help of direct growing the graph with the calibrated weight function $f(k)$ described above and the distribution of the number $x$ of arcs in the stochastic increments of the graph. The recurrent formulas for the exact calculation of the VDD in a growing graph with a stochastic increment are derived in [6]:

$$Q_g = \frac{r_g \langle f \rangle}{\langle f \rangle + mf(g)}, \qquad (27)$$

$$Q_k = \frac{r_k \langle f \rangle + mf(k-1)Q_{k-1}}{\langle f \rangle + mf(k)}, \qquad k \geq g+1, \qquad (28)$$

where $g = k_{\min}$; $r_g, r_{g+1}, \ldots, r_h$ are probabilities that the number $x$ of arcs of an increment is equal to $g$, $g+1, \ldots, h$ respectively (in our case $g = 1$, $h = 6$).

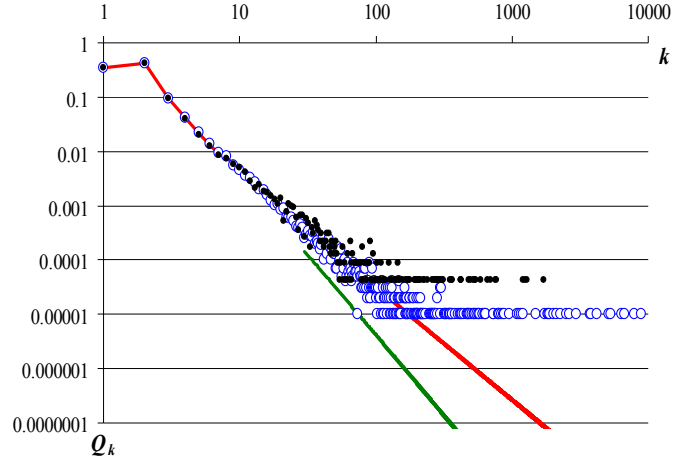

**Figure 3**. The VDD $Q_k$ of a graph with a size of 100 thousand vertices, grown in a simulation model (blue round markers) and VDD calculated by the exact formulas (27)–(28) (the upper red solid line). Black dots represent the VDD $Q_k$ in the network of autonomous systems. The green bottom solid line shows the slope angle of the VDD $Q_k$, which is realized with the traditional use of the weight function $f(k) = k$ that implements the asymptotic distribution $Q_k \sim ck^{-3}$. With this asymptotic distribution by changing the coefficient $c$ the green line can be replaced higher, but its slope angle cannot be changed.

Figure 3 clearly demonstrates that the correct choice of the weight function allows us to adequately reproduce the asymptotic VDD of the simulated networks in graph models. Due to this, the quality of solving modeling problems is significantly improved. In particular, adequate forecasting of changes in network characteristics and network processes caused by the growth of networks is provided.

Note, as an illustration, that the VDD implemented in the considered example, which has a much heavier tail than the VDD of BA graph, causes a more intensive growth of vertex degrees and a larger share of «influencers» – vertices with very high degrees. This difference is noticeable even when visually comparing the grown graphs, if they are depicted on a scale that provides the presence of several thousand vertices of the graph in one figure.

## 9. Conclusion

The article deals with the growing random graphs of $P$ class (hybrid Pennock graphs), generated using a mixture of two weight functions $f_1(k) = 1$ and $f_2(k) = k$, and random graphs of $L$ class grown using a linear weight function $f(k) = k + s$ or $f(k) = 1$.

The study proved that every hybrid graph of $P$ class has a twin – a graph of $L$ class with exactly the same final VDD, the arc endpoints degree distribution, and the edge endpoints degree distribution. A stronger statement is also proved that every $P$-graph is equivalent to its $L$-twin-graph. A formula is given that helps to define $L$-twin-graph corresponding to $P$-graph. The asymptotic VDD $Q_k$ are found for all values of weight function parameters in both classes of graphs.

It is established, that graphs with asymptotic distribution $Q_k \sim ck^{-\alpha}$, where $\alpha \in (2, 3]$, widespread in real fast growing networks, are realized by a linear weight function $f(k) = k + s$ with a displacement $s$, setting in the range $m < s \leq 0$, where $m$ is a number of arcs in a graph increment. The example of applying the obtained theoretical results given in the article confirms their validity and practical value.

The obtained results can be used in the construction of graph models of real growing networks (such as scientific collaboration networks [8], social networks [9], data networks [10], etc.). These results can be applied in the study of networks and network processes both by analytical [11–13] and imitational [14, 15] methods. This can significantly improve the adequacy and efficiency of the developed methods of interaction with real networks and methods of managing network processes.

**References**
[1]   Barabasi A L and Albert R 1999 *Emergence of scaling in random networks* (Science) vol 286 pp 509–512
[2]   Barabasi A L 2009 Scale-free networks: *A decade and beyon.* (Science) vol 325 pp 412–413
[3]   Pennock D M et al 2002 *Winners don't take all: Characterizing the competition for links on the web* (Proceedings of the National Academy of Science of the United States of America, 99(8):5207–5211)
[4]   Jackson M O 2010 *Social and Economic Networks: Models and Analysis* (Stanford University, Santa Fe Institute, CIFAR)
[5]   Zadorozhnyi V N and Yudin E B 2015 *Growing network: models following nonlinear preferential attachment rule* (Physica A: Statistical Mechanics and its Applications) vol 428 pp 111–132
[6]   Zadorozhnyi V N 2011 *Random graphs with nonlinear preferential attachment rule* (Control Sciences) no. 6 pp 2–11
[7]   Zadorozhnyi V N, Yudin E B and Yudina M N 2017 *Analytical and numerical methods of calibration for preferential attachment randon graphs* (2017 International Siberian Conference on Control and Communications (SIBCON), Astana, Kazakhstan) pp 1–6
[8]   Fortunato S et al 2018 *Science of Science* (Science) vol 359 issue 6379
[9]   Deville P et al 2016 *Scaling Identity Connects Human Mobility and Social Interactions* (PNAS 113: 26 pp 7047–7052
[10]  Lo C, Cheng J and Leskovec J 2017 *Understanding Online Collection Growth Over Time: A Case Study of Pinterest* (ACM International Conference on World Wide Web (WWW))
[11]  Karrer B and Newman M E J 2011 *Competing epidemics on complex networks* (Phys. Rev. E 84, 036106)
[12]  Boccaletti S at al 2006 *Complex networks: Structure and dynamics* (Physics Reports 424) pp 175–308
[13]  Zan Y 2018 *DSIR double-rumors spreading model in complex networks* (Chaos, Solitons and Fractals) vol 110 pp 191–202
[14]  Witten G and Poulter G 2007 *Simulations of infections diseases on networks* (Computers in Biology and Medicine) vol 37 pp. 195–205
[15]  Aksyonov K at al 2017 *Architecture of the Multi-agent Resource Conversion Processes Extended with Agent Coalitions* (IEEE International Symposium on Robotics and Intelligent Sensors, IRIS 2016 Hosei University, Tokyo, Japan, 17-20 December 2016, Code 134518, Procedia Computer Science 105) pp 221–226